\documentclass[11pt]{article}

\usepackage{epsfig,graphicx}
\usepackage{amsfonts, amssymb, latexsym} 

\def\R{\mathbb{R}}

\def\R{\mathbb{R}}

\def\0{{\bf 0}}

\date{}

\begin{document}

\title{Priority and correctness  on ``From Schoenberg to Pick-Nevanlinna: Toward a complete picture of the variogram class" by Porcu and Schilling (2011)} 

\author{
CHUNSHENG MA \\
 {\small{\em Department of Mathematics, Statistics, and Physics,  Wichita State University, 
      Wichita, KS 67260-0033, USA.}} \\
  {\small{\em E-mail: chunsheng.ma@wichita.edu}} }

\maketitle
\date{}

It is not the purpose of this correspondence to complain about that six out of seven theorems listed in Porcu and Schilling \cite{Porcu2011} belong to others, but  not to the authors  themselves. This  is simply  a call for the priority,  originality, and correctness.

\begin{itemize}
\item[(i)]
Theorem 2  listed in Porcu and Schilling \cite{Porcu2011} is not correct, for which a counterexample is
$\phi (x) = \exp ( - |x|)-1,  x \in \R$.

\item[(ii)] Lemma 9 in  Porcu and Schilling \cite{Porcu2011} is not correct, for which a counterexample is
  $$ \gamma (\xi) = \left\{ 
                                  \begin{array}{ll}
                                  0,   ~   &   ~ \xi  = {\bf 0}, \\
                                  1,   ~   &   ~ \xi  \neq {\bf 0}.
                                  \end{array} \right. $$ 

\item[(iii)]  The proof of Theorem 8 on pages 448-449 of Porcu and Schilling \cite{Porcu2011},  which is Theorem 3 (i) of Ma \cite{Ma2007Pam},  is not correct. To see this, recall that $a_1$ and $a_2$ are two given positive constants in the theorem, and that, for
the ``only if" part, 
 $(1-e^{-a_1 \gamma (\xi)})((1-e^{-a_2 \gamma (\xi)})$ is assumed to be a variogram. 
  First of all, from these assumptions, there is no guarantee to claim that
     $$ \xi  \longrightarrow  \frac{(1-e^{-a_1  \gamma (\xi)})}{a_1  } \cdot
         \frac{(1-e^{-a_2  \gamma (\xi)})}{a_2 }  $$
  is a variogram   for all $a_1, a_2 > 0$.     
  Based on  the assumption that $\gamma (\xi), \xi \in \R^d,$ is a homogeneous function, one can only infer  that
  $$ \xi  \longrightarrow  \frac{(1-e^{-a_1 \lambda \gamma (\xi)})}{a_1 \lambda } \cdot
         \frac{(1-e^{-a_2 \lambda \gamma (\xi)})}{a_2 \lambda}  $$
   is a variogram   for all $\lambda  > 0$. 
   
   Second, it is doubtful to apply  a L{\' e}vy-Khinchine representation for this case, only under the  assumption that $\gamma (\xi), \xi \in \R^d,$ is a homogeneous function.

  \item[(iv)]  Section 4 of  Porcu and Schilling \cite{Porcu2011} is not new,  and Lemma 17 of  Porcu and Schilling \cite{Porcu2011} is a special case of Theorem 3 of Ma \cite{Ma2003Spl}, where
   no continuous assumption was made.

\end{itemize}

\end{document}